\documentclass[12pt]{amsart}
\usepackage{amsthm,amssymb,amsmath,amscd}
\usepackage{mdwlist}
\usepackage[noend]{algorithmic}

\newcommand{\divides}{\mid}
\newcommand{\ndivides}{\nmid}

\newtheorem{theorem}{Theorem}[section]

\newtheorem{proposition}[theorem]{Proposition}
\theoremstyle{definition}
\newtheorem{definition}[theorem]{Definition}
\newtheorem{example}[theorem]{Example}
\newtheorem{remark}[theorem]{Remark}

\title{Prime injections and quasipolarities}

\author{Octavio A. Agust\'{\i}n-Aquino}
\address{Universidad de la Ca\~nada\\
Carretera Teotitl\'an-Nanahuatipan S/N\\
Teotitl\'an de Flores Mag\'on, Oaxaca\\
M\'exico}
\email{octavioalberto@unca.edu.mx}

\begin{document}
\begin{abstract}
Let $p$ be a prime number. Consider the injection
\[
 \iota:\mathbb{Z}/n\mathbb{Z}\to\mathbb{Z}/pn\mathbb{Z}:x\mapsto px,
 \]
and the elements $e^{u}.v:=(u,v)\in \mathbb{Z}/n\mathbb{Z}\rtimes \mathbb{Z}/n\mathbb{Z}^{\times}$ and
$e^{w}.r:=(w,r)\in \mathbb{Z}_{p n}\rtimes \mathbb{Z}_{p n}^{\times}$. Suppose $e^{u}.v\in
\mathbb{Z}/n\mathbb{Z}\rtimes \mathbb{Z}/n\mathbb{Z}^{\times}$ is seen as an automorphism of $\mathbb{Z}/n\mathbb{Z}$
by $e^{u}.v(x)=vx+u$; then $e^{u}.v$ is a \emph{quasipolarity} if it is an involution without fixed points.
In this brief note give an explicit formula for the number of quasipolarities of $\mathbb{Z}/n\mathbb{Z}$ in
terms of the prime decomposition of $n$, and we prove sufficient conditions such that $(e^{w}.r)\circ \iota
=\iota\circ (e^{u}.v)$, where $e^{w}.r$ and $e^{u}.v$ are quasipolarities.
\end{abstract}
\subjclass[2000]{11A05, 11A07}
\keywords{Quasipolarity, injection, unitary divisor}
\maketitle

\section{Some preliminaries}

Before we can state and prove the result of this paper, it is convenient to explain certain conventions, notions and notations taken
from Guerino Mazzola's monograph on mathematical musicology \cite{gM02}. First of all, the ring (or $\mathbb{Z}$-module)
$\mathbb{Z}/n\mathbb{Z}$ is a good model of the $n$-tone equally tempered scale modulo octaves \cite[Chapter 6, Section 4]{gM02}. Now consider the group
\[
 \overrightarrow{GL}(\mathbb{Z}/n\mathbb{Z}):=(\mathbb{Z}/n\mathbb{Z})\rtimes (\mathbb{Z}/n\mathbb{Z})^{\times}
\]
and let us denote an element $(u,v)\in \overrightarrow{GL}(\mathbb{Z}/n\mathbb{Z})$
with $e^{u}.v$. We have an action of $\overrightarrow{GL}(\mathbb{Z}/n\mathbb{Z})$
on $\mathbb{Z}/n\mathbb{Z}$ defined in the following way:
\[
e^{u}.v(x) = vx+u.
\]

These notations are meant to exhibit the importance of the actions of affine groups on musical objects. In particular,
the exponential notation $e^{u}$ was chosen because the composition of two translations (or transpositions,
musically speaking) is $e^{u}\circ e^{v} = e^{u+v}$. The
linear part $v$ is also musically meaningful: the best example is perhaps $v=-1$, which corresponds to the inversion of intervals
and melodies (but some more are possible).

The aforementioned action of $\overrightarrow{GL}(\mathbb{Z}/n\mathbb{Z})$
on $\mathbb{Z}/n\mathbb{Z}$ extends naturally to an action on the powerset of $\mathbb{Z}/2k\mathbb{Z}$ in a pointwise manner. A \emph{marked strong dichotomy}
is a subset $D\subseteq \mathbb{Z}/2k\mathbb{Z}$ such that there is a unique $\pi=e^{u}.v$ that satisfies
\[
\pi(D)=(\mathbb{Z}/2k\mathbb{Z})\setminus D.
\]

Marked strong dichotomies are important abstractions for the mathematical theory of counterpoint as conceived by Mazzola, for they generalize the notion
of consonance and dissonance in the standard $12$-tone and microtonal equal tuning. See \cite[part VII]{gM02}, \cite{oA09} and \cite{oA11} for further details.

\section{Polarities and quasipolarities}

The unique element $\pi$ that switches between a marked strong dichotomy $D$ and its complement is called
the \emph{polarity}. This nomenclature was chosen by Mazzola because:
\begin{quote}
[T]he traditional consonance/dissonance concept is not a polar one, since intervals are more or less consonant, for example in Euler theory and Helmholtz theory. But in musical theory, they are strictly separated into one or another category. So what is more or less in the acoustical or number[...] theories is now compressed into two bags, or poles. That is the reason. (Guerino Mazzola, personal communication, September 25th, 2013).
\end{quote}

Observe that if we
regard $\pi$ as a automorphism of $\mathbb{Z}/2k\mathbb{Z}$, we have $\pi^{2}=\mathrm{Id}_{\mathbb{Z}/2k\mathbb{Z}}$
and it has no fixed points. Thus, any $\pi=e^{u}.v$ with these two properties is called a \emph{quasipolarity}.

\begin{remark}
From the definition, it is easy to see that the notion of quasipolarity makes sense for $\mathbb{Z}/n\mathbb{Z}$
only when $n$ is even: being an involution, its cycles (regarding it as a permutation) have cardinality at most $2$, and none
can actually be of cardinality $1$, because otherwise it would have a fixed point.
\end{remark}

Let $\omega(2k)$ be the number of distinct prime factors of $2k$, with $k\geq 1$, and
\[
2k=2^{\alpha}\prod_{i=1}^{\omega(2k)-1}p_{i}^{\alpha_{i}}
\]
its prime decomposition. Suppose $2k=2^{\alpha}ab$ with $a$ coprime with $b$. It is known \cite[p. 191]{iV54} that the solutions
of $v^{2}\equiv 1 \pmod{2k}$ are given by the simultaneous solutions of the pair of congruences
\begin{equation}\label{E:Cong}
x\equiv 1 \pmod{2a},\quad x\equiv -1 \pmod{2b}
\end{equation}
if $\alpha = 1,2$, or the two pairs
\begin{equation}\label{E:Cong2}
\begin{gathered}
x\equiv 1 \pmod{2^{\alpha-1}a},\quad  x\equiv -1 \pmod{2^b},\\
x\equiv 1 \pmod{2a},\quad x\equiv -1 \pmod{2^{\alpha-1}b},
\end{gathered}
\end{equation}
if $\alpha > 2$. Moreover, in \cite[Theorem 3.1]{CDM256} it is proved that the number of affine 
parts $e^{u}$ available for an involution $v$ such that $e^{u}.v$ is a quasipolarity is
$\frac{2k}{\gcd(v-1,2k)}$, whenever
\begin{equation}\label{E:CondQP}
 2\frac{2k}{\gcd(v+1,2k)}=\gcd(v-1,2k).
\end{equation}
 
With this information, we can easily compute some values of the number $Q(2k)$ of quasipolarities of $\mathbb{Z}/2k\mathbb{Z}$
for small values of $2k$, see Table \ref{T:NQP}.
\begin{table}
\label{T:NQP}
\begin{tabular}{|c|c|c|}
\hline
$2k$ & $Q(2k)$ & List\\
\hline
$2$ & $1$ & $e^{1}.1$\\
$4$ & $3$ & $e^{2}.1,e^{\{1,3\}}.3$\\
$6$ & $4$ & $e^{3}.1,e^{\{1,3,5\}}.5$\\
$8$ & $5$ & $e^{4}.1,e^{\{1,3,5,7\}}.7$\\
$10$ & $6$ & $e^{5}.1,e^{\{1,3,5,7,9\}}.5$\\
$12$ & $12$ & $e^{6}.1,e^{\{2,6,10\}}.5,e^{\{3,9\}}.7,e^{\{1,3,5,7,9,11\}}.11$\\
$14$ & $8$ & $e^{7}.1,e^{\{1,3,5,7,9,11,13\}}.13$\\
$16$ & $9$ & $e^{8}.1,e^{\{1,3,5,7,9,11,13,15\}}.13$\\
$18$ & $10$ & $e^{9}.1,e^{\{1,3,5,7,9,11,13,15,17\}}.17$\\
$20$ & $18$ & $e^{10}.1,e^{\{2,6,10,14,18\}}.9,e^{\{5,15\}}.11,e^{2\mathbb{Z}+1}.19$\\
$22$ & $12$ & $e^{11}.1,e^{2\mathbb{Z}+1}.21$\\
$24$ & $20$ & $e^{12}.1,e^{\{3,9,15,21\}}.7,e^{\{4,12,20\}}.17,e^{2\mathbb{Z}+1}.23$\\
\hline
\end{tabular}
\caption{The number $Q(2k)$ of quasipolarities in $\mathbb{Z}/2k\mathbb{Z}$ for $1\leq k\leq 12$, and their explicit enumeration.}
\end{table}

After querying the OEIS for the sequence $Q(2k)$, I discovered the following interesting fact \cite{oeisA034448}, that was also
implied by a conjecture of the referee of this paper. We need a definition first.

\begin{definition}
A divisor $d$ of $n$ is said to be \emph{unitary} if $\gcd(d,\frac{n}{d})=1$, and we write $d||n$.
\end{definition}

\begin{proposition}
Let $\sigma^{*}_{1}(n):=\sum_{d||n}d$ be the sum of the unitary divisors of $n$. We have
\[
Q(2k) = \sigma^{*}_{1}(k).
\]
\end{proposition}

If the prime decomposition of $k$ is known, we can calculate $\sigma^{*}_{1}(k)$ in a more direct way \cite{HS79}:
\[
\sigma^{*}_{1}(k) = \prod_{p^{u}||k}(1+p^{u}).
\]

\begin{proof}
Suppose first that $\alpha=1$. We have that $v$ satisfies \eqref{E:Cong},
\[
v \equiv 1 \pmod{2a}, \quad v\equiv -1 \pmod{2b},
\]
thus $2a\divides (v-1)$ and $2b\divides (v+1)$. Now
\[
\gcd(v-1,2k) = 2a\quad\text{and}\quad \gcd(v+1,2k)=2b
\]
because $2k$ is divisible by $2$ only once. This means that
\[
 2\frac{2k}{\gcd(v+1,2k)} = 2\frac{2ab}{2a}=2b=\gcd(v+1,2k)
\]
which implies there are exactly
\[
\frac{2k}{\gcd(v-1,2k)}= \frac{2ab}{2a} = b
\]
quasipolarities of the form $e^{u}.v$. From this it is evident that each involution $v$ is in correspondence
with a unitary divisor $b$ of $k$, thus
\[
\sum_{v^{2}\equiv 1\pmod{2k}} \frac{2k}{\gcd(v-1,2k)} = \sum_{b||k} b = \sigma_{1}^{*}(k).
\]

If $\alpha=2$ the same reasoning works \emph{mutatis mutandis} if $\gcd(v-1,2k) = 4a$. Otherwise, $\gcd(v-1,2k)=2a$,
which means that $v-1=2q$ with $q$ odd. Then $v+1=2(q+1)=4q'$, hence $\gcd(v+1,2k)=4b$, which rehabilitates the
argument because $2b$ is a unitary divisor of $k$.

The case $\alpha\geq 3$ is slightly more difficult. The symmetry of the systems of congruences \eqref{E:Cong2} enable
us to suppose without loss of generality that $\gcd(v-1,2k)=2^{\beta}a$ with $\beta\geq \alpha-1$. If $\gcd(v-1,2k) = 2^{\alpha}a$, then necessarily
$\gcd(v+1,2k)=2b$, for $v\equiv -1 \pmod{2b}$ and $\gcd(v-1,v+1)=2$, and therefore the proof goes as before.
If $\gcd(v-1,2k) = 2^{\alpha-1}a$, then $v$ does not define quasipolarities, but $w=v+k$ does. Indeed,
\[
w^{2}=(v+k)^{2} = v^{2}+2vk+k^{2} =v^{2}+2kv+2^{\alpha}kab \equiv v^{2} \equiv 1 \pmod{2k},
\]
and since $v-1=2^{\alpha-1}aq$ for some odd $q$ then
\[
w-1=2^{\alpha-1}aq+k = 2^{\alpha-1}aq+2^{\alpha-1}ab = 2^{\alpha-1}a(q+b)
\]
where $q+b$ is even (being the sum of two odd integers). Thus $\gcd(w-1,2k)=2^{\alpha}a$ and $\gcd(w+1,2k)=2b$,
which yields the summand for $b$ anew.
\end{proof}

\section{Quasipolarities and injections}

In \cite[Ch. 4]{oA11} it was shown that, whenever
\begin{enumerate}
\item there is a marked strong dichotomy $D$ in $\mathbb{Z}/2k\mathbb{Z}$ with polarity $\pi$ and,
\item there is a quasipolarity in $\pi'\in\overrightarrow{GL}(\mathbb{Z}_{4k})$ such that $\iota\circ \pi = \pi'\circ \iota$ (where $\iota:\mathbb{Z}/2k\mathbb{Z}\to\mathbb{Z}_{4k}$ is the injection given by $\iota(x)=2x$),
\end{enumerate}
then $\pi'$ is the polarity of the marked strong dichotomy $2D=\{2x:x\in D\}$ in $\mathbb{Z}_{4k}$. This means that the notion of consonance (and dissonance)
can be meaningfully extended from a $2k$-tone equal tuning to a $4k$-equal tuning; in other words, it is possible to ``lift'' the original
consonances of a $2k$-tone scale to the $4k$-tone scale, in the sense that the lifted polarity restricts to original polarity in the $2k$-tone scale.

If we were to generalize this result for injections of the form
\begin{align*}
 \iota:\mathbb{Z}/n\mathbb{Z}&\to\mathbb{Z}/{p n}\mathbb{Z},\\
 x&\mapsto px,
\end{align*}
were $p$ is a prime number, we would need first to find the conditions such that $\iota \circ p = p'\circ \iota$,
which is the aim of this section.


Suppose $e^{u}.v$ is a quasipolarity. As we said before, in \cite{CDM256} it is proved that $u$ can
be taken as $\frac{2k}{\gcd(v+1,n)}$. Our interest is to find the conditions
such that there exists a quasipolarity $e^{w}. r:\mathbb{Z}/pn\mathbb{Z}\to\mathbb{Z}/pn\mathbb{Z}$ which renders
commutative the following square:
\begin{equation}\label{E:CuadComm}
\begin{CD}
\mathbb{Z}/n\mathbb{Z} @>\iota>> \mathbb{Z}_{p n}\\
@Ve^{u}.vVV @VVe^{w}.rV\\
\mathbb{Z}/n\mathbb{Z} @>>\iota> \mathbb{Z}_{p n}.
\end{CD}
\end{equation}

We have the following result.
\begin{theorem}
Let $n$ be an even number, $v$ an involution modulo $n$, $k=\frac{v^{2}-1}{n}$ and
$u=\frac{n}{\gcd(v+1,n)}$. If either $\gcd(p,2v)\divides k$ or $p\ndivides n$,
then there exists $t$ such that $r=v+nt$ is an involution modulo $p n$. In that
case and if $\frac{r^{2}-1}{p n}$ is even, then $e^{w}\circ r$ is a quasipolarity with $w=p u$ and the diagram \eqref{E:CuadComm} 
commutes.
\end{theorem}

\begin{proof}We begin noting that for
the square \eqref{E:CuadComm} to commute it is necessary and sufficient that
\[
 w \equiv p u\pmod{p n}, \quad p r \equiv p v \pmod{p n}.
\]

The second congruence is equivalent to $p(r-v) = p n t$ for some
integer $t$. Hence $r-v=nt$ and
\[
 r = v + nt.
\]

Let $v$ is an involution in $\mathbb{Z}/n\mathbb{Z}$. We want $r$ to be an
involution. We see that
\begin{align*}
 (v+nt)^2 &= v^{2}+2vnt+n^2t^2\\
 &= 1+kn+2vnt+n^2t^2\\
 &= 1 + (k+2vt+nt^2)n,
\end{align*}
so for $(v+nt)$ to be an involution is necessary and sufficient that
$p \divides (k+2vt+nt^2)$. In other words, $t$ is the solution
of the quadratic congruence
\begin{equation}\label{E:QuadCong}
 nt^{2}+2vt+k \equiv 0 \pmod{p}.
\end{equation}

We distinguish two cases. If $p \divides n$, then it is enough
to solve for $t$ the linear congruence
\[
 2vt \equiv -k \pmod{p}.
\]

Such a congruence is solvable if and only if $\gcd(p,2v)\divides -k$.
Note that if $p=2$, this condition simply means that $k$ must a multiple of $2$.

If $p\ndivides n$, the quadratic congruence is unavoidable.
Fortunately, $\gcd(2n,p)=1$
so, in order to solve it, we rewrite \eqref{E:QuadCong} to obtain
\[
(2nt+2v)^{2}\equiv 4v^{2}-4nk \pmod{p}
\]
which reduces to
\[
 (nt+v)^{2} \equiv v^{2}-nk \equiv 1+nk-nk \equiv 1\pmod{p}.
\]

Since $1$ is always a quadratic residue, we
deduce that $t = n^{-1}(\pm1-v)$ where $n^{-1}$
is the inverse of $n$ modulo $p$.

Suppose now that we have found a $t$ such that $r=v+tn$ is an involution. For there exists $w$ such that
$e^{w}\circ r$ is a quasipolarity, it is necessary and sufficient
to check that
\begin{equation}\label{E:C1}
 2\frac{p n}{\gcd(v+nt+1,p n)}=\gcd(v+nt-1,p n).
\end{equation}

If \eqref{E:C1} is true, we can choose
\begin{equation}\label{E:AffPart}
 w = \frac{p n}{\gcd(v+nt+1,p n)}.
\end{equation}

Let us begin. Note that
\[
 \gcd\left(\frac{v+nt-1}{2},\frac{v+nt+1}{2}\right)=1
\]
and
\[
 \gcd(v+nt\pm 1,p n) = 2\gcd\left(\frac{v+nt\pm 1}{2},p\frac{n}{2}\right)
\]
thus
\begin{align*}
 \gcd(v+nt+1,p n)\gcd(v+nt-1,p n) &= 4\gcd\left(\frac{(v+nt)^{2}-1}{4},p\frac{n}{2}\right)\\
 &=4\gcd\left(\frac{1+k'p n-1}{4},p \frac{n}{2}\right)\\
 &=2\gcd\left(\frac{k'p n}{2},p n\right)\\
 &= 2p\gcd\left(k'\frac{n}{2},n\right).
\end{align*}

Observing that
\[
 \gcd\left(k'\frac{n}{2},n\right)= n
\]
holds if and only if $2\divides k'$, it follows that \eqref{E:C1} holds if
and only if $2\divides k'$, where $k'=\frac{r^{2}-1}{p n}$.

To finish, we show that once $e^{w}\circ r$ is an involution such that
$rp \equiv vp\pmod{p n}$ and $w$ is given by \eqref{E:AffPart}, it
is true that $w$ equals $p u$ and thus the diagram \eqref{E:CuadComm} commutes. In case $p\ndivides (v+tn+1)$, then
\[
 \gcd(v+tn+1,p n) = \gcd(v+tn+1,n ) = \gcd(v+1,n)
\]
which means that
\begin{equation}\label{E:Prop}
 w = \frac{p n}{\gcd(v+tn+1,p n)} = \frac{p n}{\gcd(v+1,n)} = p u.
\end{equation}

Assume now the alternative case $p\divides (v+tn+1)$. Then any
divisor $d$ of $\frac{v+tn+1}{p}$ and $n$ is also a divisor of $v+1$,
because $v+1$ is a linear combination of them:
\[
 v+1 = p\frac{v+tn+1}{p}-tn.
\]

It follows that
\[
\gcd\left(\frac{v+tn+1}{p},n\right)\divides \gcd(v+1,n)
\]
or, equivalently,
\[
\gcd(v+1,n) = \gcd\left(\frac{v+tn+1}{p},n\right)p = \gcd\left(v+tn+1,p n\right).
\]

The missing factor is $p$ because any common factor $d$ of $v+1$ and $n$ divides
the linear combination $(v+1)+tn$ and also $\frac{(v+1)+tn}{p}$, as long as
$\gcd(d, p)=1$ or $d=p^{\lambda-1}$, where $p^{\lambda}$ is the
greatest power of $p$ that divides both $v+tn+1$ and $n$. In conclusion,
equation \eqref{E:Prop} is true again, and the proof concludes.
\end{proof}

\begin{example}
The affine map $e^{2}\circ 5:\mathbb{Z}_{12}\to\mathbb{Z}_{12}$ is a
quasipolarity. Let $p = 2$ and $k=\frac{5^{2}-1}{12}=2$. Since $2\divides k$,
there exists a $t$ such that $5+12t$ is an involution in $\mathbb{Z}_{24}$.
Using the proof of the theorem, $t$ is the solution of
\[
0 \equiv 2\cdot v t \equiv  -k \equiv - 2 \equiv 0 \pmod{2},
\]
thus $t$ can be chosen arbitrarily. If we choose $t=0$, $k'=\frac{5^{2}-1}{24}=1$ is not even. If
$t=1$, then $r=5+12=17$ and $k'=\frac{17^{2}-1}{24}=12$ is
even and $w=p u = 2\cdot 2 = 4$. Hence $e^{4}\circ 17:\mathbb{Z}_{24}\to\mathbb{Z}_{24}$
is a quasipolarity such that \eqref{E:CuadComm} commutes.

If now we take $p=5$, we have $5\ndivides 12$, so $t=3(\pm 1 - 5)\bmod 5 = \pm 3$. If
we choose $t=3$, we get $r=5+2\cdot 12= 29$ and it is such that $\frac{r^{2}-1}{60}=14$
is even, so $e^{10}\circ 29:\mathbb{Z}_{60}\to\mathbb{Z}_{60}$ satisfies \eqref{E:CuadComm}.
\end{example}

\section*{Acknowledgments}

I sincerely thank the referee's many observations that improved my exposition and results. I also greatly thank
Guerino B. Mazzola and Shalosh B. Ekhad for answering my questions about polarities and quasipolarities patiently and quickly.



\bibliographystyle{amsplain}
\bibliography{tesis}


\end{document}